%% file: ModelStrOperads_v3.tex
\documentclass[11pt,a4paper,english]{article}
\usepackage[latin1]{inputenc}
\usepackage[T1]{fontenc}
\usepackage{amsmath,amssymb,amscd}
\usepackage{babel,epsfig}
\usepackage{mathrsfs}
\usepackage{pifont}
\usepackage[all]{xy}

\newcommand{\lite}[1]{{\scriptstyle #1}}
\newcommand{\co}{\colon\thinspace}    
\newcommand{\Com}{\mathbf{Com}}
\newcommand{\Ass}{\mathbf{Ass}}

\newcommand{\Hopf}{\operatorname{Hopf}}
\newcommand{\Form}[1]{{}_{#1}\!\operatorname{Form}}
\newcommand{\Mod}[1]{{}_{#1}\!\operatorname{Mod}}
\newcommand{\Alg}[1]{{}_{#1}\!\operatorname{Alg}}
\newcommand{\Oper}{\operatorname{Oper}}
\newcommand{\Coll}{\operatorname{Coll}}

\newcommand{\op}{\operatorname{op}}
\newcommand{\E}{\mathcal{E}}

\newcommand{\V}{\mathcal{V}}
\newcommand{\R}{\mathbb{R}}
\newcommand{\Z}{\mathbb{Z}}

\newcommand{\Sp}{\mathcal{S}p}

\newcommand{\hocolim}{\operatornamewithlimits{hocolim}}

\newcommand{\colim}{\operatornamewithlimits{colim}}

\newtheorem{Def}{Definition}[section]

\newtheorem{Thm}[Def]{Theorem}
\newtheorem{Prop}[Def]{Proposition}

\newtheorem{Cor}[Def]{Corollary}

\newtheorem{Rem}[Def]{Remark}

\newenvironment{proof}{\mbox{}\newline\textbf{Proof:} }{\nopagebreak\hfill$\square$\\}

\title{Model structure on operads in orthogonal spectra}
\author{Tore August Kro\\
\texttt{toreak@math.ntnu.no}}

\begin{document}

\maketitle

\section*{Abstract}

We generalize Berger and Moerdijk's results on axiomatic homotopy
theory for operads to the setting of enriched symmetric monoidal
model categories, and show how this theory applies to orthogonal
spectra. In particular, we provide a symmetric fibrant replacement
functor for the positive stable model structure.

\input{introduction}

\input{body_v3}

\bibliographystyle{amsalpha}
\bibliography{MinBibliografi,Ekstrabib}

\end{document}

%% file: introduction.tex
\section{Introduction}

Operads (in topological spaces) were introduced in order to describe
algebraic structures where the constraints are relaxed up to a
system of homotopies. The definition of operads generalizes to any
symmetric monoidal category. This raises the question about
axiomatic homotopy theory for operads, given that the base category
has a monoidal model structure. This question has answers when the
base category is simplicial sets by Rezk~\cite{Rezk:96}, complexes
of a module over a ring by Hinich~\cite{Hinich:97}, a cofibrantly
generated symmetric monoidal model category by
Spitzweck~\cite{Spitzweck:01}, $k$-spaces by Vogt~\cite{Vogt:03},
and a closed symmetric monoidal model category by Berger and
Moerdijk~\cite{BergerMoerdijk:03}.

The aim of this article is to provide Quillen model structures on
operads and their algebras when the base category is some symmetric
monoidal category of spectra, for instance orthogonal spectra,
see~\cite[example~4.4]{MandellMaySchwedeShipley:01}. An argument of
Lewis~\cite{Lewis:91} shows that no symmetric monoidal model
category of spectra can simultaneously have a cofibrant unit and a
symmetric monoidal fibrant replacement functor. Unfortunately,
Berger and Moerdijk's constructions require both these properties of
the base category.

We resolve this problem by considering an enriched symmetric
monoidal model category. This generalizes the hint given
in~\cite[example~4.6.4]{BergerMoerdijk:03}. We follow the strategy
and proofs of the paper~\cite{BergerMoerdijk:03} closely, and
observe that there is no need to assume that the unit of the base
category is cofibrant, as long as the category we enrich in has a
cofibrant unit. We advise the reader to keep a copy of Berger and
Moerdijk's article at hand while reading section~\ref{modcat}.

The category of orthogonal spectra, with the positive stable model
structure, see~\cite[section~14]{MandellMaySchwedeShipley:01},
satisfies a priori nearly all of the requirements of the previous
section. The only missing piece is a symmetric monoidal fibrant
replacement functor. We show in section~\ref{caseSpO} that the
second-most naive guess for a fibrant replacement functor actually
is symmetric. The analogous functor in symmetric spectra is not a
fibrant replacement.

\textbf{Acknowledgements.} The author would like to thank K. Hess
for her encouragement and detailed comments in this research. Also
thanks to M. Ching for email discussions, to H. Fausk for comments,
and to I. Moerdijk for comments and for explaining some details of
their article.

%% file: body_v3.tex

\section{Model structures on operads and their left
modules}\label{modcat}

We assume that the reader is familiar with the basic notions of
model categories, see for example~\cite{DwyerSpalinski:95}
or~\cite{Hirschhorn:03}.

For the basic definitions of enriched categories,
consult~\cite{Kelly:05}. We recall a few facts: Throughout $\V$
denotes a closed symmetric monoidal category with product $\otimes$,
unit $I$, and \emph{internal hom} $[-,-]$. Let $\E$ be a symmetric
monoidal category with product $\wedge$, unit $S$, enriched with hom
objects $\E(-,-)$ in $\V$, having tensor
$\odot\co\V\times\E\rightarrow\E$, and cotensor written by
exponentiation. For $A,B$ in $\E$ and $Y$ in $\V$ there are natural
isomorphisms
\[
[Y,\E(A,B)]\cong\E(Y\odot A,B)\cong \E(A,B^Y)\quad.
\]

We will assume that both $\V$ and the underlying category of $\E$
are monoidal model categories, see~\cite{Hovey:98}
and~\cite{SchwedeShipley:00}, in particular the
\emph{pushout-product axiom} holds. Moreover, we need an axiom
relating the model structures on $\V$ and $\E$. Thus we assume the
\begin{itemize}
\item[]\textbf{Pullback-cotensor axiom:} If $p\co A\twoheadrightarrow B$
is an $\E$-fibration and $i\co X\hookrightarrow Y$ is a
$\V$-cofibration, then $p^{\square i}\co A^Y\rightarrow
B^Y\times_{B^X}A^X$ is an $\E$-fibration, and moreover $p^{\square
i}$ is trivial if either $p$ or $i$ is trivial.
\end{itemize}
By adjunction this axiom have two equivalent reformulations, one of
them similar to Quillen's axiom \textbf{SM7}. Moreover, our axiom
implies that $\E$ is a model enrichment by $\V$, in the sense
of~\cite[section~3.1]{Dugger:06}.

Let $\Hopf(\V)$ be the category of \emph{commutative Hopf objects}
in $\V$, see~\cite[section~1]{BergerMoerdijk:03}. Observe that any
Abelian monoid $M$ naturally gives rise to a commutative Hopf object
$I[M]$ whose underlying object in $\V$ is $\coprod_M I$. Consider
$\Z/2$ multiplicatively. If the folding map $I[\Z/2]\rightarrow I$
can be factored in $\Hopf(\V)$ as $I[\Z/2]\hookrightarrow
H\xrightarrow{\simeq} I$, where the underlying maps in $\V$ are a
cofibration and a weak equivalence respectively, then we say that
$\V$ \emph{admits a commutative Hopf interval}.

For a finite group $G$, let $\E^G$ denote the category of objects in
$\E$ with a right $G$-action and $G$-equivariant maps. A
\emph{collection} in $\E$ is a sequence of objects $A(n)$, $n\geq0$
in $\E$, such that $A(n)$ has a right action of the symmetric group
$\Sigma_n$. This category $\Coll(\E)$ equals the product
$\prod_{n=0}^\infty\E^{\Sigma_n}$. Assuming that $\E$ is cofibrantly
generated, there is a model structure on collections where
$A\rightarrow B$ is a weak equivalence (resp. fibration) if each
$A(n)\rightarrow B(n)$ is a non-equivariant weak equivalence (resp.
fibration) in $\E$. The subcategory $\widetilde{\Coll}(\E)$ (resp.
$\Coll_+(\E)$) of \emph{reduced collections} (resp. \emph{positive
collections}) consists of those $A$ such that $A(0)=S$ (resp.
$A(0)=\emptyset$).

An \emph{operad in $\E$} is a collection $\mathcal{P}$ together with
a unit $S\rightarrow\mathcal{P}(1)$ and structure maps
\[
\mathcal{P}(k)\wedge\mathcal{P}(n_1)\wedge\cdots\wedge\mathcal{P}(n_k)\rightarrow\mathcal{P}(n_1+\cdots+n_k)
\]
satisfying certain conditions, see~\cite{May:72}. Alternatively, one
can define this category $\Oper(\E)$ as the monoids for Smirnov's
non-commutative monoidal product on collections in $\E$,
see~\cite{Smirnov:82}
or~\cite[section~I.1.8]{MarklShniderStasheff:02}. We denote this
product by $\circ$, and will define and study it more thoroughly
later in this note, see definition~\ref{Def:circprod}. The unit for
$\circ$ is the collection $\mathcal{S}$ with $\mathcal{S}(1)=S$ and
$\mathcal{S}(n)=\emptyset$ for $n\neq1$. An operad $\mathcal{P}$ is
called \emph{positive} if $\mathcal{P}(0)=\emptyset$, or is
\emph{reduced} if $\mathcal{P}(0)=S$. Denote these categories
$\Oper_+(\E)$ and $\widetilde{\Oper}(\E)$ respectively. Let $\Ass$
and $\Com$ denote the operads for associative and commutative
monoids respectively. Their $n$-ary parts are given by
$\Ass(n)=S[\Sigma_n]$ and $\Com(n)=S$. Observe that the category of
reduced operads is the subcategory of $\Oper(\E)/\Com$ consisting of
$\alpha\co\mathcal{P}\rightarrow\Com$ with $\alpha(0)$ being the
identity of $S$. An operad $\mathcal{P}$ is called
\emph{$\Sigma$-split} if $\mathcal{P}$ is a retract of
$\mathcal{P}\wedge\Ass$.

For an arbitrary operad $\mathcal{P}$ we define categories
$\Mod{\mathcal{P}}$, $\Alg{\mathcal{P}}$, and
$\Form{\mathcal{P}}^d$. A \emph{left $\mathcal{P}$-module} is a
collection $M$ together with a left action $\mathcal{P}\circ
M\rightarrow M$. A \emph{$\mathcal{P}$-algebra} $A$ is a left
$\mathcal{P}$-module concentrated in arity $0$, i.e. $A(n)=*$ for
$n>0$. Explicitly, we have structure maps $\mathcal{P}(n)\wedge
A^{\wedge n}\rightarrow A$. More generally, we define a \emph{$d$'th
order $\mathcal{P}$-form} to be a left $\mathcal{P}$-module
truncated above arity $d$, i.e. $M(n)=*$ for $n>d$. A
\emph{$\mathcal{P}$-coalgebra} is an object $B$ of $\E$ together
with structure maps $B\wedge\mathcal{P}(n)\rightarrow B^{\wedge n}$
satisfying conditions dual to those of a $\mathcal{P}$-algebra.

Let $\mathcal{D}$ be one of the categories $\Oper_+(\E)$,
$\widetilde{\Oper}(\E)$, $\Mod{\mathcal{P}}$,
$\Form{\mathcal{P}}^d$, or $\Alg{\mathcal{P}}$. In all cases we have
forgetful functors to $\Coll(\E)$. We say that $\mathcal{D}$
\emph{admits a transferred model structure} if there is a model
structure on $\mathcal{D}$ where $A\rightarrow B$ is a weak
equivalence (resp. fibration) if and only if the underlying map in
$\Coll(\E)$ is a weak equivalence (resp. fibration).

Generalizing the main results of \cite{BergerMoerdijk:03} we have:

\begin{Thm}
Assume that the unit $I\in\V$ is cofibrant, $\V$ admits a
commutative Hopf interval, $\E$ is cofibrantly generated, and $\E/S$
(resp. $\E$) admits a symmetric monoidal fibrant replacement
functor. The category of reduced operads in $\E$ (resp. positive
operads in $\E$) then admits a transferred model structure.
\end{Thm}

\begin{Thm}\label{thm:modstronmod}
Assume that the unit $I\in\V$ is cofibrant, $\E$ is cofibrantly
generated, and $\E$ has a symmetric fibrant replacement functor. Let
$\mathcal{P}$ be an operad in $\E$ and $\mathcal{Q}$ an operad in
$\V$. If there exists an operad map
$j\co\mathcal{P}\rightarrow\mathcal{Q}\odot\mathcal{P}$ and an
interval in $\V$ with a $\mathcal{Q}$-coalgebra structure, then
$\Mod{\mathcal{P}}$, $\Form{\mathcal{P}}^d$ and $\Alg{\mathcal{P}}$
admit transferred model structures.
\end{Thm}

\begin{Cor}
Assume that the unit $I\in\V$ is cofibrant, $\E$ is cofibrantly
generated, and $\E$ has a symmetric fibrant replacement functor. If
there exists an interval in $\V$ with a coassociative
comultiplication, then for all $\Sigma$-split operads $\mathcal{P}$
the categories $\Mod{\mathcal{P}}$, $\Form{\mathcal{P}}^d$ and
$\Alg{\mathcal{P}}$ admit transferred model structures.
\end{Cor}

\begin{Cor}
Assume that the unit $I\in\V$ is cofibrant, $\E$ is cofibrantly
generated, and $\E$ has a symmetric fibrant replacement functor. If
there exists an interval in $\V$ with a coassociative and
cocommutative comultiplication, then for all operads $\mathcal{P}$
the categories $\Mod{\mathcal{P}}$, $\Form{\mathcal{P}}^d$ and
$\Alg{\mathcal{P}}$ admit transferred model structures.
\end{Cor}

\begin{proof}
We prove all four results simultaneously and follow Berger and
Moerdijk closely in their approach. Hence, we will only outline the
arguments, to the extent it becomes obvious that everything they do
also work in our enriched setting.

Since $\Alg{\mathcal{P}}$ and $\Form{\mathcal{P}}^d$ are truncations
of $\Mod{\mathcal{P}}$, we will not mention them again in this
proof, i.e. the details are exactly as for left modules. Moreover,
the two corollaries follow from theorem~\ref{thm:modstronmod} by
taking $\mathcal{Q}=\Ass$ and $\Com$ respectively.

To put model structures on $\Oper_+(\E)$, $\widetilde{\Oper}(\E)$
and $\Mod{\mathcal{P}}$, we consider free-forgetful adjunctions
\begin{align*}
\Coll_+(\E)&\rightleftarrows\Oper_+(\E),\\
\widetilde{\Coll}(\E/S)&\rightleftarrows\widetilde{\Oper}(\E),\quad\text{and}\\
\Coll(\E)&\rightleftarrows\Mod{\mathcal{P}}.
\end{align*}
Using the transfer principle and Quillen's path-object argument, as
explained in~\cite[sections~2.5 and~2.6]{BergerMoerdijk:03}, we have
to check that $\Oper_+(\E)$, $\widetilde{\Oper}(\E)$ and
$\Mod{\mathcal{P}}$ have small colimits and finite limits, the free
functors preserve small objects, $\Oper_+(\E)$,
$\widetilde{\Oper}(\E)$ and $\Mod{\mathcal{P}}$ have fibrant
replacement functors, and $\Oper_+(\E)$, $\widetilde{\Oper}(\E)$ and
$\Mod{\mathcal{P}}$ have functorial path-objects for fibrant
objects. See also~\cite[theorem~11.3.2]{Hirschhorn:03}.

The functorial fibrant replacement functors of $\Oper_+(\E)$,
$\widetilde{\Oper}(\E)$ and $\Mod{\mathcal{P}}$ are defined
aritywise using the symmetric fibrant replacement functors of $\E$,
$\E/S$ and $\E$ respectively. To get the functorial path-objects we
use convolution pairings
\begin{align*}
\Hopf(\V)^{\op}\times\Oper_+(\E) &\rightarrow\Oper_+(\E),\\
\Hopf(\V)^{\op}\times\widetilde{\Oper}(\E) &\rightarrow\widetilde{\Oper}(\E),\quad\text{and}\\
\operatorname{Coalg}_\mathcal{Q}^{\op}\times\Mod{\mathcal{P}}
&\rightarrow\Mod{\mathcal{Q}\odot\mathcal{P}}\quad.
\end{align*}

To construct the first pairing, observe that each commutative Hopf
object $H$ defines a cooperad $TH$ with $TH(n)=H^{\otimes n}$. We
use the unreduced convolution pairing $\mathcal{P}^{TH}$ given by
$\mathcal{P}^{TH}(n)=\mathcal{P}(n)^{TH(n)}$. Now let $H$ be a
commutative Hopf interval in $\V$. As
in~\cite[theorem~3.1]{BergerMoerdijk:03} we get a functorial path
object
\[
\mathcal{P}=\mathcal{P}^{TI}\xrightarrow{\simeq}\mathcal{P}^{TH}
\twoheadrightarrow \mathcal{P}^{TI[\Z/2]}\twoheadrightarrow
\mathcal{P}\times\mathcal{P}
\]
for positive operads.

For the second pairing we are given $H\in\Hopf(\V)$ and
$\mathcal{P}\rightarrow\Com$. Observe that the counit
$\epsilon:H\rightarrow I$ is a map of commutative Hopf objects.
Thus, we can define the reduced convolution pairing
$\widetilde{\mathcal{P}^{TH}}$ as the pullback of
\[
\Com\xrightarrow{\epsilon^*}\Com^{TH}\leftarrow\mathcal{P}^{TH}\quad.
\]
Let $\mathcal{P}$ be a fibrant reduced operad in $\E$, and let
$I[\Z/2]\hookrightarrow H\xrightarrow{\simeq}I$ be a commutative
Hopf interval in $\V$. By convolution we get
\[
\mathcal{P}=\widetilde{\mathcal{P}^{TI}}\xrightarrow{\simeq}\widetilde{\mathcal{P}^{TH}}
\twoheadrightarrow
\widetilde{\mathcal{P}^{TI[\Z/2]}}\twoheadrightarrow
\mathcal{P}\times_\Com\mathcal{P}\quad.
\]
$H$ and $I$ are cofibrant, so the first map is a weak equivalence by
Ken Brown's lemma and the cotensor-pullback axiom. The middle map is
a fibration by the cotensor-pullback axiom. While the last map is
for $n\geq1$ the projection
$\widetilde{\mathcal{P}^{I[\Z/2]}}(n)=\mathcal{P}(n)^{\times_S
2^n}\rightarrow \mathcal{P}(n)\times_S\mathcal{P}(n)$ onto the first
and last factor, and hence a fibration. This yields a functorial
path-object for fibrant $\mathcal{P}$.

The last convolution pairing, $M^B$, between a coalgebra $B$ under
an operad $\mathcal{Q}$ in $\V$ and a left $\mathcal{P}$-module $M$,
is defined by the formula $M^B(n)=M(n)^B$. Here, the left
$\mathcal{Q}\odot\mathcal{P}$-module structure map
\[
(\mathcal{Q}\odot\mathcal{P})(k)\wedge M^B(n_1)\wedge\cdots\wedge
M^B(n_k)\rightarrow M^B(n)
\]
is given as the adjoint of the composition
\begin{equation*}
\begin{gathered}
B\odot(\mathcal{Q}\odot\mathcal{P})(k)\wedge
M^B(n_1)\wedge\cdots\wedge M^B(n_k)\\
\begin{aligned}
\quad\quad\quad\quad\quad\quad\quad&\cong
\left(B\otimes\mathcal{Q}(k)\right)\odot\left(\mathcal{P}(k)\wedge
M(n_1)^B\wedge\cdots\wedge M(n_k)^B\right)\\
&\rightarrow B^{\otimes k}\odot\left(\mathcal{P}(k)\wedge
M(n_1)^B\wedge\cdots\wedge M(n_k)^B\right)\\
&\cong \mathcal{P}(k)\wedge
(B\odot M(n_1)^B)\wedge\cdots\wedge (B\odot M(n_k)^B)\\
&\rightarrow \mathcal{P}(k)\wedge M(n_1)\wedge\cdots\wedge
M(n_k)\rightarrow M(n)
\end{aligned}
\end{gathered}\quad.
\end{equation*}
From here on the path-object argument works exactly as in Berger and
Moerdijk's proof of~\cite[theorem~4.1]{BergerMoerdijk:03}.
\end{proof}

The initial positive operad is $\mathcal{S}$, the unit for the
$\circ$-product on collections, while the initial reduced operad
$\widetilde{\mathcal{S}}$ is given by $\widetilde{\mathcal{S}}(n)=S$
for $n=0,1$ and $\widetilde{\mathcal{S}}(n)=\emptyset$ otherwise. If
the unit $S$ of $\E$ is not cofibrant, we have to be extra careful
with cofibrancy of operads. Therefore, we call a reduced operad
$\mathcal{P}$ \emph{$\Sigma$-cofibrant} if the unique map
$\widetilde{\mathcal{S}}\rightarrow\mathcal{P}$ is a cofibration in
the model structure on collections. Similarly, we also call a
positive operad $\mathcal{P}$ \emph{$\Sigma$-cofibrant} if the
unique map $\mathcal{S}\rightarrow\mathcal{P}$ is a cofibration of
collections. An map of operads $\mathcal{P}\rightarrow\mathcal{Q}$
is called a \emph{$\Sigma$-cofibration} if the underlying map of
collections is a cofibration. Observe that our notions of
$\Sigma$-cofibrant differs from the definition found
in~\cite[section~4]{BergerMoerdijk:03}, but agrees with the
definition in~\cite[section~2.4]{BergerMoerdijk:06}. However, all
notions coincide if $S$ is cofibrant in $\E$.

\begin{Prop}
Any cofibrant reduced (resp. positive) operad is $\Sigma$-cofibrant.
\end{Prop}

\begin{proof}
We consider the case of reduced operads first. Since the initial
reduced operad, $\widetilde{S}$, is $\Sigma$-cofibrant, it is enough
to show that $\Sigma$-cofibrant reduced operads are closed under
cellular extensions. We will now contemplate the difference between
reduced and unreduced operads. Let
$F\co\Coll(\E)\rightarrow\Oper(\E)$ be the free operad functor, and
let
$\tilde{F}\co\widetilde{\Coll}(\E/S)\rightarrow\widetilde{\Oper}(\E)$
be the free reduced operad functor. Given a reduced collection $A$
in $\E/S$, we observe that $\tilde{F}A(n)=FA(n)$ for $n>0$, while
$\tilde{F}A(0)=S$.

\cite[corollary~5.2]{BergerMoerdijk:03} says that for any
cofibration $A\hookrightarrow B$ of collections and any map of
operads $FA\rightarrow\mathcal{P}$ the induced map
$\mathcal{P}\rightarrow\mathcal{P}\cup_{FA}FB$ is a
$\Sigma$-cofibration. So, if $\mathcal{P}$ is a $\Sigma$-cofibrant
reduced operad, $A\hookrightarrow B$ a cofibration in
$\widetilde{\Coll}(\E/S)$, and $u\co A\rightarrow U(\mathcal{P})$ an
arbitrary map, then the only difference between
$\mathcal{P}\cup_{\tilde{F}A}\tilde{F}B$ and
$\mathcal{P}\cup_{FA}FB$ lies in arity $0$. Hence, the map
$\mathcal{P}\rightarrow \mathcal{P}\cup_{\tilde{F}A}\tilde{F}B$ is a
$\Sigma$-cofibration.

Next, we consider the case of positive operads. Observe that the
free positive operad functor is the restriction of
$F:\Coll(\E)\rightarrow\Oper(\E)$ to the subcategory of positive
collections. Hence, \cite[corollary~5.2]{BergerMoerdijk:03}
immediately applies, and yields the result.
\end{proof}

We now turn towards Smirnov's product on collections,
see~\cite{Smirnov:82}. We begin with the definition of the
$\circ$-product, and proceed by proving a two technical results,
namely the propositions~\ref{prop:circpropertyI}
and~\ref{prop:circpropertyII}. Berger and Moerdijk prove
technicalities of similar flavor
in~\cite[section~2.5]{BergerMoerdijk:06}.

\begin{Def}\label{Def:circprod}
Let $r_1,\ldots,r_k$ be non-negative integers that sum up to $n$.
Abbreviate $X(r_1)\wedge\cdots\wedge X(r_k)$ by $X(r_*)$, and
$\Sigma_{r_1}\times\cdots\times\Sigma_{r_k}$ by $\Sigma(r_*)$. We
get a $\Sigma(r_*)$-action on $X(r_*)$. Now induce up to a
$\Sigma_n$-equivariant object
$\operatorname{Ind}^{\Sigma_n}_{\Sigma(r_*)}X(r_*)$. On the disjoint
union over all partitions of length $k$ that sum to $n$,
\[
\coprod_{r_*}\operatorname{Ind}^{\Sigma_n}_{\Sigma(r_*)}X(r_*)\quad,
\]
we have a $\Sigma_k$-action by permuting factors in $X(r_*)$ and
blocks in $\Sigma_n$. Now define
\[
(A\circ X)(n)=\coprod_{k=0}^\infty
A(k)\wedge_{\Sigma_k}\left(\coprod_{r_*}\operatorname{Ind}^{\Sigma_n}_{\Sigma(r_*)}X(r_*)\right)\quad.
\]
\end{Def}

We will now derive a few properties of this product, but before that
let us introduce a piece of terminology coming from Goodwillie's
calculus of functors: We call the diagram
\[
\xymatrix@=12pt{
A \ar[r]\ar[d] & B \ar[d] \\
C \ar[r] & D }
\]
\emph{a cofibration square} if the three maps $A\rightarrow B$,
$A\rightarrow C$ and $B\cup_A C\rightarrow D$ are cofibrations.

\begin{Prop}\label{prop:circpropertyI}
Let $A$ be a cofibrant in the model category of reduced collection
under $\widetilde{\mathcal{S}}$ (resp. positive collections under
$\mathcal{S}$), and let $X\hookrightarrow Y$ be a cofibration
between cofibrant collections. Then $A\circ X\rightarrow A\circ Y$
is also a cofibration.
\end{Prop}

\begin{proof}
Fix $n$ and $k$. $\Sigma_k$ acts on partitions $r_*$ of length $k$
that sum up to $n$. Fix also a representative $r_*$ for an orbit of
this action. Let $\operatorname{Aut}(r_*)$ be the permutations in
$\Sigma_k$ that acts trivially on $r_*$. It is enough to show that
\[
A(k)\wedge_{\operatorname{Aut}(r_*)}\left(\operatorname{Ind}^{\Sigma_n}_{\Sigma(r_*)}X(r_*)\right)
\rightarrow
A(k)\wedge_{\operatorname{Aut}(r_*)}\left(\operatorname{Ind}^{\Sigma_n}_{\Sigma(r_*)}Y(r_*)\right)
\]
is a $\Sigma_n$-equivariant cofibration. Clearly, $X(r_*)\rightarrow
Y(r_*)$ is a $\Sigma(r_*)$-equivariant cofibration, and
$\operatorname{Ind}_{\Sigma(r_*)}^{\Sigma_n}$ preserves
(equivariant) cofibrations. Observe that the map above is
tautologically an equivariant cofibration if $A$ is the initial
reduced operad $\widetilde{\mathcal{S}}$ (resp. the initial positive
operad $\mathcal{S}$). In general, we may assume that $A(k)$ is a
cellular $\E^{\Sigma_k}$-object relative to
$\widetilde{\mathcal{S}}(k)$ (resp. rel $\mathcal{S}(k)$). Hence
$A(k)$ is a (transfinite) sequential colimit where each step
$A(k)_\alpha\hookrightarrow A(k)_{\alpha+1}$ is formed by gluing a
generating cofibration $\Sigma_k\times
\partial_\alpha\hookrightarrow\Sigma_k\times D_\alpha$. By the pushout product axiom
for $\E$ we have a cofibration square
\[
\xymatrix{
\coprod_{\Sigma_k/{\operatorname{Aut}(r_*)}}\partial_\alpha\wedge
\left(\operatorname{Ind}^{\Sigma_n}_{\Sigma(r_*)}X(r_*)\right)
\ar[r] \ar[d] &
\coprod_{\Sigma_k/{\operatorname{Aut}(r_*)}}D_\alpha\wedge
\left(\operatorname{Ind}^{\Sigma_n}_{\Sigma(r_*)}X(r_*)\right) \ar[d] \\
\coprod_{\Sigma_k/{\operatorname{Aut}(r_*)}}\partial_\alpha\wedge
\left(\operatorname{Ind}^{\Sigma_n}_{\Sigma(r_*)}Y(r_*)\right)
\ar[r] & \coprod_{\Sigma_k/{\operatorname{Aut}(r_*)}}D_\alpha\wedge
\left(\operatorname{Ind}^{\Sigma_n}_{\Sigma(r_*)}Y(r_*)\right)}
\]
Hence, also
\[
\xymatrix{ A(k)_\alpha\wedge_{\operatorname{Aut}(r_*)}
\left(\operatorname{Ind}^{\Sigma_n}_{\Sigma(r_*)}X(r_*)\right)
\ar[r] \ar[d] & A(k)_{\alpha+1}\wedge_{\operatorname{Aut}(r_*)}
\left(\operatorname{Ind}^{\Sigma_n}_{\Sigma(r_*)}X(r_*)\right) \ar[d] \\
A(k)_\alpha\wedge_{\operatorname{Aut}(r_*)}
\left(\operatorname{Ind}^{\Sigma_n}_{\Sigma(r_*)}Y(r_*)\right)
\ar[r] &A(k)_{\alpha+1}\wedge_{\operatorname{Aut}(r_*)}
\left(\operatorname{Ind}^{\Sigma_n}_{\Sigma(r_*)}Y(r_*)\right)}
\]
is a cofibration square. The conclusion follows.
\end{proof}

\begin{Prop}\label{prop:circpropertyII}
Let $A\xrightarrow{\simeq}B$ be a weak equivalence between cofibrant
reduced collections under $\widetilde{\mathcal{S}}$ (resp. cofibrant
positive collections under $\mathcal{S}$), and let $X$ be a
cofibrant collection. Then $A\circ X \rightarrow B\circ X$ is also a
weak equivalence.
\end{Prop}

\begin{proof}
By Ken Brown's lemma, it is enough to consider the case where
$A\rightarrow B$ is an acyclic cofibration between reduced (resp.
positive) collections. Fix $n$, $k$ and $r_*$ as above. We may
assume that $B(k)$ is cellular relative to $A(k)$. Hence, we write
$B(k)$ as a (transfinite) sequential colimit starting with
$B(k)_0=A(k)$, and such that each step $B(k)_{\alpha}\rightarrow
B(k)_{\alpha+1}$ is the pushout along a generating acyclic
cofibration, $\Sigma_k\times \partial_\alpha\hookrightarrow
\Sigma_k\times D_\alpha$. We now get a pushout diagram
\[
\xymatrix{
\coprod_{\Sigma_k/\operatorname{Aut}(r_*)}\partial_\alpha\wedge
\left(\operatorname{Ind}^{\Sigma_n}_{\Sigma(r_*)}X(r_*)\right)
\ar[r] \ar[d] &
\coprod_{\Sigma_k/\operatorname{Aut}(r_*)}D_\alpha\wedge
\left(\operatorname{Ind}^{\Sigma_n}_{\Sigma(r_*)}X(r_*)\right) \ar[d] \\
B(k)_\alpha\wedge_{\operatorname{Aut}(r_*)}
\left(\operatorname{Ind}^{\Sigma_n}_{\Sigma(r_*)}X(r_*)\right)
\ar[r] &B(k)_{\alpha+1}\wedge_{\operatorname{Aut}(r_*)}
\left(\operatorname{Ind}^{\Sigma_n}_{\Sigma(r_*)}X(r_*)\right)}
\]
where the top map is an acyclic cofibration by the pushout-product
axiom. Hence the bottom map also is an acyclic cofibration.
\end{proof}

The free functor
$F_\mathcal{P}\co\Coll(\E)\rightarrow\Mod{\mathcal{P}}$ is given by
the $\circ$-product, i.e. $F_\mathcal{P}(X)=\mathcal{P}\circ X$.
This extends the Schur functor defining the free
$\mathcal{P}$-algebra.

\begin{Thm}
Under the assumptions of theorem~\ref{thm:modstronmod}, assume
additionally that $\E$ is left proper, and that the domains of the
generating cofibrations are cofibrant. If
$\phi\co\mathcal{P}\rightarrow\mathcal{Q}$ is a weak equivalence
between $\Sigma$-cofibrant reduced (resp. positive) operads, then
the base-change adjunctions
$\Mod{\mathcal{P}}\rightleftarrows\Mod{\mathcal{Q}}$,
$\Alg{\mathcal{P}}\rightleftarrows\Alg{\mathcal{Q}}$, and
$\Form{\mathcal{P}}^d\rightleftarrows\Form{\mathcal{Q}}^d$ are
Quillen equivalences.
\end{Thm}

\begin{proof}
We prove this for left modules, the other cases are similar. The
categories $\Mod{\mathcal{P}}$ and $\Mod{\mathcal{Q}}$ both carry
transferred model structures, hence the base-change adjunction is a
Quillen pair by inspection. Since the forgetful functor $\phi^*$
reflects weak equivalences, it is enough to show that the unit of
the adjunction, $M\rightarrow\phi^*\phi_! M$, is a weak equivalence
for each cellular left $\mathcal{P}$-module $M$. Cellular means that
$M$ is a (transfinite) sequential colimit starting from the initial
left $\mathcal{P}$-module, and where $M_{\alpha+1}$ is the pushout
of $M_\alpha\leftarrow \mathcal{P}\circ
\partial_\alpha\hookrightarrow \mathcal{P}\circ D_\alpha$, for a
generating cofibration $\partial_\alpha\rightarrow D_\alpha$ in
$\Coll(\E)$. Observe $\phi^*\phi_!M$ inherits a similar description,
i.e. $\phi^*\phi_!M_{\alpha+1}$ is the pushout of
$\phi^*\phi_!M_\alpha\leftarrow \phi^*(\mathcal{Q}\circ
\partial_\alpha)\rightarrow \phi^*(\mathcal{Q}\circ D_\alpha)$. Recall from
proposition~\ref{prop:circpropertyII} that $\mathcal{P}\circ
\partial_\alpha\rightarrow \mathcal{Q}\circ \partial_\alpha$ and $\mathcal{P}\circ
D_\alpha\rightarrow \mathcal{Q}\circ D_\alpha$ are weak
equivalences, while $\mathcal{P}\circ \partial_\alpha\rightarrow
\mathcal{P}\circ D_\alpha$ and $\mathcal{Q}\circ
\partial_\alpha\rightarrow \mathcal{Q}\circ D_\alpha$ are
$\Sigma$-cofibrations by proposition~\ref{prop:circpropertyI}. Thus
inductively, all $M_\alpha\rightarrow \phi^*\phi_!M_\alpha$ are weak
equivalences, and the conclusion follows.
\end{proof}

\section{The case of orthogonal spectra}\label{caseSpO}

It is convenient to replace the category of all topological spaces
by compactly generated spaces (= weak Hausdorff $k$-spaces,
see~\cite{McCord:69}). We define the reduced homotopy colimit of a
sequence
$X_1\xrightarrow{f_1}X_2\xrightarrow{f_2}X_3\rightarrow\cdots$ of
based spaces as the reduced mapping telescope. $\hocolim_n X_n$ has
the topology of the union $\bigcup_{n=0}^\infty F_n$, where the
$n$'th space of the filtration is
\[
F_n=\left(X_1\wedge I_+\right) \cup_{f_1} \left(X_2\wedge I_+\right)
\cup_{f_2} \left(X_3\wedge I_+\right) \cup \cdots \cup
\left(X_{n-1}\wedge I_+\right) \cup_{f_{n-1}} X_n\quad.
\]
Since each $X_n$ is compactly generated, all base points $*\in X_n$
are closed, so any compact subset of $\hocolim_n X_n$ is contained
in some $F_k$, see~\cite[lemma~9.3]{Steenrod:67}. Using the
projections $F_n\rightarrow X_n$, it is easy to prove that we have
natural group isomorphisms
\[
\colim_n\pi_q X_n\xrightarrow{\cong}\pi_q\hocolim_n X_n\quad.
\]

An \emph{orthogonal spectrum} $X$ consists of a based
$O(V)$-equivariant space $X(V)$ for every finite-dimensional real
inner product space $V$ together with $O(V)\times O(W)$-equivariant
suspension maps $\sigma\co S^V\wedge X(W)\rightarrow X(V\oplus W)$
satisfying the obvious coherence condition,
see~\cite[example~4.4]{MandellMaySchwedeShipley:01}. Fixing
orthogonal spectra $X$ and $Y$ consider pairs $(Z,\mu)$, where $Z$
is an orthogonal spectrum and $\mu$ is a family of maps $\mu(V,W)\co
X(V)\wedge Y(W)\rightarrow Z(V\oplus W)$ such that $\mu(V,W)$ is
$O(V)\times O(W)$-equivariant and the following diagram commutes for
all $U$, $V$ and $W$:
\begin{equation}\tag{*}
\xymatrix@C=2.2cm{ & {\scriptstyle X(U\oplus V)\wedge Y(W)} \ar[dr]^{\mu(U\oplus V,W)} & \\
{\scriptstyle S^U\wedge X(V) \wedge Y(W)} \ar[ur]^{\sigma_X\wedge 1}
\ar[r]_{1\wedge\mu(V,W)} \ar[d]_{\text{twist}} & {\scriptstyle
S^U\wedge Z(V\oplus W)} \ar[r]_{\sigma_Z} & {\scriptstyle Z(U\oplus V\oplus W)} \\
{\scriptstyle X(V)\wedge S^U\wedge Y(W)} \ar[r]^{1\wedge\sigma_Y} &
{\scriptstyle X(V)\wedge Y(U\oplus W)} \ar[r]^{\mu(V,U\oplus W)} &
{\scriptstyle Z(V\oplus U\oplus W)} \ar[u]_{\cong}}
\end{equation}
The \emph{smash product} $X\wedge Y$ is initial among such $Z$.
Thus, any such pair $(Z,\mu)$ determines a unique map $X\wedge
Y\xrightarrow{\mu} Z$. We denote the category of orthogonal spectra
by $\Sp^O$. The smash product $\wedge$ is symmetric monoidal with
unit the sphere spectrum $S$. Moreover, $\Sp^O$ is enriched,
tensored and cotensored over topological spaces. The \emph{stable
homotopy groups} of an orthogonal spectrum $X$ are defined for all
integers $q$ in terms of homotopy groups of topological spaces by
the formula $\pi_qX=\colim_n\pi_{q+n}X(\R^{n})$. A map $X\rightarrow
Y$ of orthogonal spectra is a \emph{weak equivalence} if it induces
isomorphisms $\pi_qX\cong\pi_qY$ of all $q$. This definition of weak
equivalence fits into a model structure on $\Sp^O$:

\begin{Thm}[{\cite[section~14]{MandellMaySchwedeShipley:01}}]
There is a model structure on $\Sp^O$ called the \emph{positive
stable model structure}, where the weak equivalences are as above.
The model structure is cofibrantly generated, left and right proper,
and satisfies the pullback-cotensor and pushout-product axioms.
Furthermore, the domains of the generating cofibrations are
cofibrant. The fibrations are characterized as the maps
$E\rightarrow B$ such that for all $V$ of positive dimension
$E(V)\rightarrow B(V)$ is a Serre fibration and the diagram
\[
\xymatrix{ E(V) \ar[r] \ar[d] & \Omega E(V\oplus \R) \ar[d] \\ B(V)
\ar[r] & \Omega B(V\oplus \R)}
\]
is homotopy pullback.
\end{Thm}

\begin{Def}
Abbreviate $\R^n\otimes V$ by $\lite{n}V$. Define the functor $T$ by
the formula $TX(V)=\hocolim_n\Omega^{nV}X(\lite{(n+1)}V)$.
\end{Def}

Let $\mathscr{K}$ be the category of pairs $(V,W)$, where the
morphisms $(V,W)\rightarrow(V',W')$ consists of injective linear
isometries $i\co V\rightarrow V'$ and $j\co W\rightarrow W'$
together with a linear isometric isomorphism $\alpha\co V'-i(V)\cong
W'-j(W)$. If $V,W,V',W'$ are oriented, then $V'-i(V)$ and $W'-j(W)$
inherit orientations, and we call a morphism in $\mathscr{K}$
\emph{positively oriented} if $\alpha$ is orientation preserving.
Thus we have a category $\mathscr{K}_+$ of pairs of oriented inner
product spaces and positively oriented morphisms. It is easily seen
that the space $\mathscr{K}_+\big( (V,W),(V',W')\big)$ is connected
whenever $\dim V<\dim V'$.

Now observe that each orthogonal spectrum $X$ gives rise to a
continuous functor from $\mathscr{K}$ to spaces by sending $(V,W)$
to $\Omega^VX(W)$. Here the induced maps $\Omega^VX(W)\rightarrow
\Omega^{V'}X(W')$ come from applying $\Omega^V(-)$ to the adjoint of
$\sigma\co S^U\wedge X(W)\rightarrow X(U\oplus W)$ with $U=V'-i(V)$.

Abbreviate $\Omega^{\lite{n}V}X(\lite{(n+1)}V)$ by $\Omega^{\otimes
n}X(V)$. The morphism
\[
(\lite{n}W,\lite{(n+1)}W)\rightarrow (\lite{(n+1)}V\oplus
\lite{n}W,\lite{(n+1)}V\oplus \lite{(n+1)}W)
\]
induces a map
\[
\Omega^{\lite{n}W}X(\lite{(n+1)}W)\rightarrow\Omega^V\Omega^{\lite{n}(V\oplus
W)}X(\lite{(n+1)}(V\oplus W))
\]
whose adjoint makes $\Omega^{\otimes n}X$ into an orthogonal
spectrum. The morphisms
$(\lite{n}V,\lite{(n+1)}V)\rightarrow(\lite{n}V\oplus V,
\lite{(n+1)}V\oplus V)$ induce the structure maps of the homotopy
colimit:
\[
\Omega^{\otimes n}X(V)
=\Omega^{\lite{n}V}X(\lite{(n+1)}V)\rightarrow\Omega^{\lite{(n+1)}V}X(\lite{(n+2)}V)
=\Omega^{\otimes(n+1)}X(V)\quad.
\]
Moreover, these maps commute with the suspension for the orthogonal
spectra $\Omega^{\otimes n}X$ and $\Omega^{\otimes(n+1)}X$.

Fix orthogonal spectra $X$ and $Y$, let $k=\max(m,n)$, and define
$\mu_{m,n}(V,W)$ as the composition
\[
\begin{gathered}
\Omega^{\lite{m}V}X(\lite{(m+1)}V)\wedge \Omega^{\lite{n}W}Y(\lite{(n+1)}W)\\
\begin{aligned}
\quad\quad\quad\quad\quad&\rightarrow\Omega^{\lite{k}V}X(\lite{(k+1)}V)\wedge
\Omega^{\lite{k}W}Y(\lite{(k+1)}W)\\
&\rightarrow\Omega^{\lite{k}V\oplus \lite{k}W}\big(X(\lite{(k+1)}V)\wedge Y(\lite{(k+1)}W)\big)\\
&\rightarrow\Omega^{\lite{k}(V\oplus W)}(X\wedge
Y)(\lite{(k+1)}(V\oplus W))\quad.
\end{aligned}
\end{gathered}
\]
This family of $\mu_{m,n}$'s satisfy (*), and hence they induce a
map $\mu_{m,n}\co\Omega^{\otimes m}X\wedge\Omega^{\otimes
n}Y\rightarrow\Omega^{\otimes k}(X\wedge Y)$. By inspection the
following four diagrams commute:
\[
\xymatrix@C=3cm{ \Omega^{\otimes n_1}X \wedge \Omega^{\otimes n_2}Y
\wedge \Omega^{\otimes n_3}Z \ar[d]^{\mu_{n_1,n_2}\wedge 1}
\ar[r]^-{1\wedge\mu_{n_2,n_3}} & \Omega^{\otimes n_1}X \wedge
\Omega^{\otimes k_2}(Y\wedge Z) \ar[d]^{\mu_{n_1,k_2}} \\
\Omega^{\otimes k_1}(X\wedge Y)\wedge \Omega^{\otimes n_3}Z
\ar[r]^{\mu_{k_1,n_3}} & \Omega^{\otimes k}(X\wedge Y\wedge Z)}
\]
where $k_1=\max(n_1,n_2)$, $k_2=\max(n_2,n_3)$, and
$k=\max(n_1,n_2,n_3)$.
\[
\xymatrix{ \Omega^{\otimes n}X \wedge \Omega^{\otimes m}Y
\ar[d]^{\text{twist}} \ar[r]^{\mu_{n,m}} & \Omega^{\otimes
k}(X\wedge Y) \ar[d]^{\text{twist}} \\ \Omega^{\otimes m}Y \wedge
\Omega^{\otimes n}X \ar[r]^{\mu_{m,n}} & \Omega^{\otimes k}(Y\wedge
X)}
\]
\[
\xymatrix{ \Omega^{\otimes n}X \wedge \Omega^{\otimes m}Y
\ar[r]^{\mu_{n,m}} \ar[d] & \Omega^{\otimes k}(X\wedge Y) \\
\Omega^{\otimes n+1}X \wedge \Omega^{\otimes m}Y
\ar[ur]_{\mu_{n+1,m}}} \quad\text{for $n<m$}
\]
\[
\xymatrix{ \Omega^{\otimes n}X \wedge \Omega^{\otimes m}Y
\ar[r]^{\mu_{n,m}} \ar[d] & \Omega^{\otimes k}(X\wedge Y) \ar[d] \\
\Omega^{\otimes n+1}X \wedge \Omega^{\otimes m}Y
\ar[r]^{\mu_{n+1,m}} & \Omega^{\otimes k+1}(X\wedge Y) }
\quad\text{for $n\geq m$.}
\]

\begin{Thm}
$T$ is a symmetric fibrant replacement functor for the positive
stable model structure on orthogonal spectra.
\end{Thm}

\begin{proof}
Clearly, $TX=\hocolim_n\Omega^{\otimes n}X$ becomes an orthogonal
spectrum, and comes with a coaugmentation $X\rightarrow TX$. By the
diagrams above the $\mu_{n,m}$'s yields a natural transformation
\[
\mu\co TX\wedge TY\rightarrow T(X\wedge Y)\quad.
\]
Moreover, $\mu$ is associative and commutative. Hence $T$ is a
symmetric monoidal functor. It remains to show that $TX$ is fibrant
for all $X$, and that the map $X\rightarrow TX$ is a weak
equivalence.

Assume $\dim V>0$. The $q$'th homotopy group of $TX(V)$ is
calculated as $\colim_n\Omega^{\lite{n}V}X(\lite{(n+1)}V)$. Now
consider the following solid diagram in $\mathscr{K}_+$:
\[
\xymatrix@=0.3cm{ \big(\lite{n}V,\lite{(n+1)}V\big) \ar[rr] \ar[dd]
\ar@{}[rd]|-{\circ} &&
\big(\lite{(2n+1)}V,\lite{(2n+2)}V\big) \ar[dd] \\ & \mbox{} \\
\big(\R\oplus \lite{n}(\R\oplus V),\lite{(n+1)}(\R\oplus V)\big)
\ar[rr] \ar@{.>}[uurr] && \big(\R\oplus \lite{(2n+1)}(\R\oplus
V),\lite{(2n+2)}(\R\oplus V)\big) \ar@{}[ul]|-{\simeq} }
\]
Choosing a linear isometry $\R\rightarrow V$, we explicitly
construct the dotted arrow such that the upper triangle commutes.
However, the lower triangle will only commute up to a homotopy. This
uses that the morphism spaces of $\mathscr{K}_+$ are connected. Thus
the induced diagram
\[
\xymatrix{ \pi_q\Omega^{\lite{n}V}X\big(\lite{(n+1)}V\big) \ar[r]
\ar[d]
 &
\pi_q\Omega^{\lite{(2n+1)}V}X\big(\lite{(2n+2)}V\big) \ar[d] \\
\pi_q\Omega\Omega^{\lite{n}(\R\oplus V)}X\big(\lite{(n+1)}(\R\oplus
V)\big) \ar[r] \ar[ur] & \pi_q\Omega\Omega^{\lite{(2n+1)}(\R\oplus
V)}X\big(\lite{(2n+2)}(\R\oplus V)\big)}
\]
commutes. This shows that $TX(V)\rightarrow\Omega TX(\R\oplus V)$ is
a weak equivalence. Hence, $TX(V)$ is fibrant in the positive stable
model structure.

A similar argument shows that the map $X\rightarrow TX$ is a weak
equivalence.
\end{proof}

\begin{Rem}
The construction of the functor $T$ can also be a carried out in the
category of symmetric spectra. However, in this case $TX$ will in
general not be a positive $\Omega$-spectrum. For a counterexample
consider $X=F_1S^1$, the free symmetric spectrum generated by a
circle $S^1$ in level $1$. The reason that the proof fails for
symmetric spectra is that the category corresponding to
$\mathscr{K}_+$ will be discrete, and hence allows no non-trivial
homotopies.

In order to construct a symmetric fibrant replacement functor for
symmetric spectra other techniques are required. Perhaps one might
improve on the small object argument.
\end{Rem}

\begin{Thm}
For an orthogonal spectrum $X$ over $S$ let $\tilde{T}X$ be the
levelwise homotopy pullback of $S\rightarrow TS\leftarrow TX$. This
defines a symmetric fibrant replacement functor $\tilde{T}$ for the
positive stable model structure on orthogonal spectra over $S$.
\end{Thm}

\begin{proof}
By construction $\tilde{T}X\rightarrow S$ is a level fibration. Fix
some positive dimensional $V$. Consider the diagram
\[
\xymatrix{ \tilde{T}X(V) \ar[r] \ar[d] & TX(V) \ar[d]
\ar[r]^-{\simeq} & \Omega TX(V\oplus\R) \ar[d] \\ S(V) \ar[r] &
TS(V) \ar[r]^-{\simeq} & \Omega TS(V\oplus \R) }\quad.
\]
The left square is homotopy pullback by definition of $\tilde{T}X$,
while the right square is homotopy pullback since the top and bottom
maps are weak equivalences. Consequently, the outer square is
homotopy pullback. Now look at the diagram
\[
\xymatrix{ \tilde{T}X(V) \ar[r] \ar[d] & \Omega\tilde{T}X(V\oplus\R)
\ar[d] \ar[r] & \Omega TX(V\oplus\R) \ar[d] \\ S(V) \ar[r] & \Omega
S(V\oplus \R) \ar[r] & \Omega TS(V\oplus \R) }\quad.
\]
The outer squares of this and the previous diagram are the same.
Since $\Omega(-)$ commutes with homotopy pullback, and by the
definition of $\tilde{T}X$, the right square is homotopy pullback.
It follows that the left square is homotopy pullback, hence we are
done.
\end{proof}

%% file: ModelStrOperads_v3.bbl
\providecommand{\bysame}{\leavevmode\hbox to3em{\hrulefill}\thinspace}
\providecommand{\MR}{\relax\ifhmode\unskip\space\fi MR }
\providecommand{\MRhref}[2]{%
  \href{http://www.ams.org/mathscinet-getitem?mr=#1}{#2}
}
\providecommand{\href}[2]{#2}
\begin{thebibliography}{MMSS01}

\bibitem[BM03]{BergerMoerdijk:03}
Clemens Berger and Ieke Moerdijk, \emph{Axiomatic homotopy theory for operads},
  Comment. Math. Helv. \textbf{78} (2003), no.~4, 805--831. \MR{MR2016697
  (2004i:18015)}

\bibitem[BM06]{BergerMoerdijk:06}
\bysame, \emph{The {B}oardman-{V}ogt resolution of operads in monoidal model
  categories}, Topology \textbf{45} (2006), no.~5, 807--849. \MR{MR2248514}

\bibitem[DS95]{DwyerSpalinski:95}
W.~G. Dwyer and J.~Spali{\'n}ski, \emph{Homotopy theories and model
  categories}, Handbook of algebraic topology, North-Holland, Amsterdam, 1995,
  pp.~73--126. \MR{MR1361887 (96h:55014)}

\bibitem[Dug06]{Dugger:06}
Daniel Dugger, \emph{Spectral enrichments of model categories}, Homology,
  Homotopy Appl. \textbf{8} (2006), no.~1, 1--30 (electronic). \MR{MR2205213
  (2006k:55034)}

\bibitem[Hin97]{Hinich:97}
Vladimir Hinich, \emph{Homological algebra of homotopy algebras}, Comm. Algebra
  \textbf{25} (1997), no.~10, 3291--3323. \MR{MR1465117 (99b:18017)}

\bibitem[Hir03]{Hirschhorn:03}
Philip~S. Hirschhorn, \emph{Model categories and their localizations},
  Mathematical Surveys and Monographs, vol.~99, American Mathematical Society,
  Providence, RI, 2003. \MR{MR1944041 (2003j:18018)}

\bibitem[Hov98]{Hovey:98}
Mark Hovey, \emph{Monoidal model categories}, arXiv:math.AT/9803002, February
  1998.

\bibitem[Kel05]{Kelly:05}
G.~M. Kelly, \emph{Basic concepts of enriched category theory}, Repr. Theory
  Appl. Categ. (2005), no.~10, vi+137 pp. (electronic), Reprint of the 1982
  original [Cambridge Univ. Press, Cambridge; MR0651714]. \MR{MR2177301}

\bibitem[Lew91]{Lewis:91}
L.~Gaunce Lewis, Jr., \emph{Is there a convenient category of spectra?}, J.
  Pure Appl. Algebra \textbf{73} (1991), no.~3, 233--246. \MR{MR1124786
  (92f:55008)}

\bibitem[May72]{May:72}
J.~P. May, \emph{The geometry of iterated loop spaces}, Springer-Verlag,
  Berlin, 1972, Lectures Notes in Mathematics, Vol. 271. \MR{54 \#8623b}

\bibitem[McC69]{McCord:69}
M.~C. McCord, \emph{Classifying spaces and infinite symmetric products}, Trans.
  Amer. Math. Soc. \textbf{146} (1969), 273--298. \MR{MR0251719 (40 \#4946)}

\bibitem[MMSS01]{MandellMaySchwedeShipley:01}
M.~A. Mandell, J.~P. May, S.~Schwede, and B.~Shipley, \emph{Model categories of
  diagram spectra}, Proc. London Math. Soc. (3) \textbf{82} (2001), no.~2,
  441--512. \MR{2001k:55025}

\bibitem[MSS02]{MarklShniderStasheff:02}
Martin Markl, Steve Shnider, and Jim Stasheff, \emph{Operads in algebra,
  topology and physics}, Mathematical Surveys and Monographs, vol.~96, American
  Mathematical Society, Providence, RI, 2002. \MR{MR1898414 (2003f:18011)}

\bibitem[Rez96]{Rezk:96}
Charles~W. Rezk, \emph{Spaces of {A}lgebra {S}tructures and {C}ohomology of
  {O}perads}, {PhD} thesis, Massachusetts {I}nstitute of {T}echnology, May
  1996.

\bibitem[Smi82]{Smirnov:82}
V.~A. Smirnov, \emph{On the cochain complex of topological spaces}, Math. USSR
  Sbornik \textbf{43} (1982), no.~1, 133--144, Translated by P.S. Landweber.
  \MR{MR618592 (83a:55016)}

\bibitem[Spi01]{Spitzweck:01}
Markus Spitzweck, \emph{Operads, algebras and modules in general model
  categories}, arXiv:math.AT/0101102, January 2001.

\bibitem[SS00]{SchwedeShipley:00}
Stefan Schwede and Brooke~E. Shipley, \emph{Algebras and modules in monoidal
  model categories}, Proc. London Math. Soc. (3) \textbf{80} (2000), no.~2,
  491--511. \MR{MR1734325 (2001c:18006)}

\bibitem[Ste67]{Steenrod:67}
N.~E. Steenrod, \emph{A convenient category of topological spaces}, Michigan
  Math. J. \textbf{14} (1967), 133--152. \MR{35 \#970}

\bibitem[Vog03]{Vogt:03}
R.~M. Vogt, \emph{Cofibrant operads and universal {$E\sb \infty$} operads},
  Topology Appl. \textbf{133} (2003), no.~1, 69--87. \MR{MR1996461
  (2004e:18014)}

\end{thebibliography}
